\documentclass[10pt]{amsart}
\usepackage{fullpage, amssymb}

\newtheorem{theorem}{Theorem}

\newtheorem{lemma}[theorem]{Lemma}
\newtheorem{corollary}[theorem]{Corollary}

\theoremstyle{definition}

\newtheorem{example}[theorem]{Example}
\newtheorem*{terminology}{Terminology}

\theoremstyle{remark}
\newtheorem{remark}[theorem]{Remark}

\numberwithin{equation}{section}
\numberwithin{theorem}{section}

\newcommand{\A}{\mathfrak{A}}

\newcommand{\B}{\mathcal{B}}
\newcommand{\C}{\mathbb{C}}

\newcommand{\h}{\mathfrak{H}}
\newcommand{\I}{\mathcal{I}}

\newcommand{\M}{\mathcal{M}}
\newcommand{\N}{\mathcal{N}}

\newcommand{\RR}{\mathcal{R}}

\newcommand{\T}{\mathbb{T}}

\newcommand{\U}{\mathcal{U}}

\begin{document}

\title[Automorphisms of ultrapowers]{Notes on automorphisms of ultrapowers of $\text{II}_1$ factors}
\author{David Sherman}
\address{Department of Mathematics\\ University of Virginia\\ P.O. Box 400137\\ Charlottesville, VA 22904}
\email{dsherman@virginia.edu}

\begin{abstract}
In functional analysis, approximative properties of an object become precise in its ultrapower.  We discuss this idea and its consequences for automorphisms of $\text{II}_1$ factors.  Here are some sample results: (1) an automorphism is approximately inner if and only if its ultrapower is $\aleph_0$-locally inner; (2) the ultrapower of an outer automorphism is always outer; (3) for unital *-homomorphisms from a separable nuclear C*-algebra into an ultrapower of a $\text{II}_1$ factor, equality of the induced traces implies unitary equivalence.  All statements are proved using operator algebraic techniques, but in the last section of the paper we indicate how the underlying principle is related to theorems of Henson's positive bounded logic.
\end{abstract}

\maketitle

\section{Introduction}

We start with some comments on the historical record.

The general ultrapower/ultraproduct construction originates in model theory, with \L o\'s's theorem in 1955 (\cite{L}) and a wave of logical applications in the 1960s.  Ultrapowers appropriate for functional analysis appeared formally around 1970 in two main flavors: a normed version for structures such as Banach spaces, and a tracial version for finite von Neumann algebras.  It is well-known that the mathematics underlying the tracial ultrapower construction was written down much earlier in Sakai's 1962 notes (\cite[Section II.7]{Sa}), although Sakai does not use ultrapower terminology.  But it seems to be less appreciated that Sakai's write-up was motivated by a 1954 article of Wright (\cite{W}).  A student of Kaplansky, Wright worked with AW*-algebras and lattices, which may explain his paper's diminished legacy.  (Currently its most recent citation on MathSciNet is from 1974.)  Nonetheless the modern reader will easily recognize Wright's descriptions of maximal ideals and the resulting quotients as the underpinnings of (the AW*-version of) the tracial ultrapower.  Thus one may justifiably say that the tracial ultrapower is older than its ``classical" set theoretic cousin.

Still in functional analysis the role of an ultrapower is simultaneously analytic and logical.
\begin{quote}
\textbf{Pattern.} An approximative property of a structure associated to a normed space corresponds to a stronger, precise version of the property in an ultrapower.  When an ultrapower has an approximative property, it automatically acquires the precise version.
\end{quote}
We do not assert this as a theorem.  Aided by appropriately restrictive definitions, it can be supported by metamathematical results, and we give some sampling of this in Section \ref{S:logic}.  There are several examples of this pattern in the literature, of which probably the best-known is that for Banach spaces $E$ and $F$, we have that $E$ is finitely representable in $F$ if and only if $E$ embeds isometrically in an ultrapower of $F$ (\cite{HM,St}, see also \cite{He} for an analyst-friendly exposition of Banach space ultraproducts).  But its implementation is not always straightforward: given one member of an ``approximative/precise property" pair, it may not be clear what the other is.

The main body of this paper concerns new examples of this pattern, proved without logical theorems or terminology.  Our primary objects are automorphisms of $\text{II}_1$ factors and their tracial ultrapowers.  We ask: If an automorphism has a certain property, what can we say about its ultrapower?  And what can we say about automorphisms of ultrapower algebras in general?  We will return to discussion of ``approximative vs. precise" at the end of this Introduction; here let us mention that one of the precise properties involves local innerness, which was introduced recently in \cite{S2007a}.  We also note that two of our results improve on conclusions of Haagerup and St\o rmer -- but they had different motivations and proved only what they needed in the context of their long article \cite{HSpinn}.  We start by reviewing the main constructions and terminology.

\smallskip

Throughout $\M$ and $\N$ will be von Neumann algebras.  We always assume $\M$ to be a $\text{II}_1$ factor, but there are no cardinality assumptions unless explicitly stated.  Any $\text{II}_1$ factor has a unique tracial state which we denote by $\tau$, avoiding subscripts when context makes the ambient algebra clear.  The $L^2$ norm on a $\text{II}_1$ factor is $\|x\|_2 = \sqrt{\tau(x^*x)}$, and it induces the strong topology on bounded subsets.  We write $\U(\N)$ for the unitary group of $\N$.

Let $\omega$ be a free ultrafilter on $\mathbb{N}$, which one may choose to think of as an element of $(\beta \mathbb{N} \setminus \mathbb{N})$.  Set $\I_\omega \subset \ell^\infty(\M)$ to be the two-sided ideal of sequences $(x_j)$ with $\|x_j\|_2 \to 0$  as $j \to \omega$.  Then $\I_\omega$ is a maximal ideal of $\ell^\infty(\M)$, and the quotient $(\ell^\infty(\M)/\I_\omega) \triangleq \M^\omega$ is a $\text{II}_1$ factor.  We call $\M^\omega$ a \textit{tracial ultrapower} of $\M$.  It is big -- even its maximal abelian *-subalgebras fail to be countably-generated (\cite[Proposition 4.3]{P1981}).  Elements of $\M^\omega$ will be denoted either by a capital letter, e.g. $X$, or by a sequence $(x_j) \in \ell^\infty(\M)$ representing the coset, following convention by omitting ``$+ \I_\omega$."  For any self-adjoint element, positive element, projection, or unitary in $\M^\omega$, we can and always do choose a representing sequence in which all terms have this same property.  (This has been proved in many places -- see \cite[Theorem 4.10]{HL} for a very general result.)  We also identify $\M$ with the subalgebra of constant sequences.  Interest in $\M^\omega$ has largely focused on the relative commutant $\M' \cap \M^\omega$, which is the algebra of $\omega$-central sequences.  See, for instance, the celebrated papers of McDuff (\cite{McD}) and Connes (\cite{C1976}).  Ge and Hadwin showed (\cite[Theorem 3.2]{GH}) that if $\M$ is countably-generated and one assumes the continuum hypothesis, the inclusion $\M \subset \M^\omega$ is independent (up to isomorphism) of the choice of $\omega$.

There are many variations of this.  Sakai actually showed that the quotient of a finite von Neumann algebra by a maximal ideal is a finite factor, so one may construct ultrapowers over larger index sets.  One may also replace $\ell^\infty(\M)$ by a direct sum of arbitrary finite factors $\{\M_j\}$, in which case the quotient by a maximal ideal is called a \textit{tracial ultraproduct} of the $\{\M_j\}$.   When $\N$ is not finite, one may still define $\I_\omega$ as the bounded sequences which go to 0 *-strongly as $n \to \omega$, but this is not an ideal of $\ell^\infty(\N)$.  Appropriate generalizations of $\M' \cap \M^\omega$ and $\M^\omega$ which are valid for arbitrary $\N$ were first defined by Connes (\cite[Section II]{C1974}) and Ocneanu (\cite[Chapter 5.1]{O}), respectively.  Over in the category of C*-algebras, ultrapowers are defined essentially as they are for Banach spaces: change $\I_\omega$ to the closed ideal of sequences which converge to 0 in norm as $n \to \omega$.  And finally, there is the model theoretic ultrapower, in which $\I_\omega$ is the algebraic ideal in $\ell^\infty(\M)$ consisting of sequences which are 0 in a neighborhood of $\omega$.  We have listed all these constructions mostly to remind the reader what we are \textit{not} doing.

Now we turn to automorphisms.  Here are some ways in which an automorphism $\theta$ of a $\text{II}_1$ factor can be ``close to inner."
\begin{enumerate}
\item We say $\theta$ is \textit{pointwise inner} if, on any self-adjoint element, it agrees with some inner automorphism.
\item We say $\theta$ is \textit{locally inner} if, on any element, it agrees with some inner automorphism.  More generally, for a cardinal $\kappa$ we say that $\theta$ is \textit{$\kappa$-locally inner} if, on any set of $\leq \kappa$ elements, it agrees with some inner automorphism.
\item We say $\theta$ is \textit{approximately inner} if, for any finite set $\{x_k\}$ and $\varepsilon > 0$, there is a unitary $u$ with
$$\max_k \|\theta(x_k) - u x_k u^*\|_2 < \varepsilon.$$
\end{enumerate}

We should immediately remark that the ``real" definition of pointwise innerness, due to Haagerup and St\o rmer (\cite[Definition 12.3]{HSeq}) and applicable to any von Neumann algebra, is that for any normal state $\varphi$, there is a unitary $u$ with $\varphi = \varphi \circ \text{Ad}(u)$.  Thus it is a predual version of local innerness (although it came first, and it uses only the positive part of the predual).  In a $\text{II}_1$ factor, this is equivalent to the definition above (\cite[Lemma 2.2]{HSpinn}), which may be thought of as ``$\frac12$-local innerness."  Now if $\kappa$ is the minimal cardinality of a generating set, and $\kappa > \kappa' > \kappa''$, we have
\begin{equation} \label{E:impl}
\text{inner} \Leftrightarrow \text{$\kappa$-locally inner} \Rightarrow \text{$\kappa'$-locally inner} \Rightarrow \text{$\kappa''$-locally inner} \Rightarrow \text{locally inner} \Rightarrow \text{pointwise inner}.
\end{equation}
Somewhat surprisingly, in a countably-generated $\text{II}_1$ factor these distinctions are meaningless: pointwise innerness already implies innerness (\cite[Proposition 12.5]{HSeq}).  But for general $\text{II}_1$ factors, none of the one-way implications in \eqref{E:impl} can be reversed (\cite[Proposition 2.2, Theorem 2.5, and Section 4.2]{S2008}), except that we have no examples to distinguish the classes of $\kappa$-locally inner automorphisms, $1\leq \kappa \leq \aleph_0$.  They are all the same if every countably-generated von Neumann algebra is singly-generated, which is the famous \textit{generator problem}.

Approximate innerness in a $\text{II}_1$ factor $\M$ says that $\theta$ belongs to the closure of the inner automorphisms in the point-strong topology.  Recall that Murray and von Neumann defined a $\text{II}_1$ factor $\M$ to have \textit{property $\Gamma$} if for any $\varepsilon > 0$ and finite set $\{x_j\} \subset \M$, there is a unitary $u$ with $\tau(u) = 0$ and $\max_j \|ux_j - x_j u \|_2 < \varepsilon$ (\cite[Definition 6.1.1]{MvN}).  Property $\Gamma$ is approximative; we discuss it further in Example \ref{T:Gamma}.  It is mentioned here because for countably-generated $\M$, the inner automorphisms are point-strong closed (and thus approximately inner implies inner) if and only if $\M$ does not have $\Gamma$ (\cite{Sa1974}, \cite[Section III]{C1974}).  Examples: free group factors do not have $\Gamma$, while the hyperfinite factor, which we denote throughout the paper as $\RR$, does.  Approximate innerness is evidently implied by $\aleph_0$-local innerness, so an affirmative answer to the generator problem would mean that approximate innerness is implied by local innerness.

\smallskip

Here are our main results.  Some of the terms are not defined precisely until later in the text.

\begin{itemize}
\item An automorphism of a $\text{II}_1$ factor is approximately inner if and only if the ultrapower of the automorphism is $\aleph_0$-locally inner (Theorem \ref{T:ai}).  Also, when an automorphism of an ultrapower is approximately inner, it is already $\aleph_0$-locally inner (Corollary \ref{T:precise}).
\item The ultrapower of an outer automorphism is always outer (Theorem \ref{T:inner}).
\item When two unital *-homomorphisms from a separable C*-algebra to an ultrapower are weakly approximately unitarily equivalent, they are already unitarily equivalent (Theorem \ref{T:ue}).  If in addition the C*-algebra is nuclear, unitary equivalence follows merely from equality of the induced traces (Corollary \ref{T:nuc}).
\item Any automorphism of an ultrapower is pointwise inner, but it need not be locally inner (Corollary \ref{T:pinn}(2)).
\end{itemize}

Now we return to ``approximative vs. precise," immediately replacing these terms with more accurate ones.

\begin{terminology}
Let $P$ and $Q$ be properties which are meaningful in some class of normed space structures (Banach spaces, tracial von Neumann algebras with automorphism, etc.), and suppose that a structure has property $P$ if and only if any ultrapower (based on a free ultrafilter of $\mathbb{N}$) has property $Q$.  Then we will say that $Q$ is the \textbf{ultrapower version} of $P$, or equivalently, $P$ is the \textbf{model version} of $Q$.
\end{terminology}

We hesitate to call this a definition, as we have not defined ``normed space structures" either in general or in the specific cases of interest.  In any event this terminology is only used to talk about the results in Sections \ref{S:aut} and \ref{S:ae}, not to state or prove them.

Our pattern, when it applies, says that $Q$ implies $P$, and that in an ultrapower they are equivalent.

\begin{remark} ${}$
\begin{enumerate}
\item The pattern does not always apply.  Suppose that $Q$ is the ultrapower version of $P$, and is strictly stronger.  Then [not $Q$] is the ultrapower version of [not $P$], and is strictly weaker.  Loosely speaking, the pattern applies to properties which assert existence, as opposed to non-existence (``positive" properties, see Section \ref{S:logic}).
\item Any property can be chosen to play either role, but there is the caveat that properties which are not equivalent in general may be equivalent in ultrapower structures.  Indeed, when the pattern applies we can always take the ultrapower version of a property $P$ to be $P$ itself.  Of course this is not interesting; we would like to name the ultrapower version as the property which is strongest in general.
\item These concepts would not change if we used any countably incomplete ultrafilter.  (An ultrafilter is \textit{countably incomplete} if it has a countable collection of members with empty intersection.)  But they would trivialize if we allowed all ultrafilters, as any model is an ultrapower of itself.  See \cite[Section 6]{GH} for more on the distinct behavior of ultrapowers based on ultrafilters which are not countably incomplete.
\item This terminology is related to, but distinct from, some concepts in Banach space theory.  Our ultrapower versions need not be ``super" in the sense that they pass to all closed subspaces of any ultrapower.  Our model versions need not be ``local" in the sense that all ultrapowers have a common property (but this would be true if we only allowed countably incomplete ultrafilters), or even ``local" in the sense that they involve approximation or finite subsets.
\item It might be more accurate to use ``ultraroot version" in place of ``model version," as of course ultrapowers are still models.  We take the perspective that ultrapowers are somehow auxiliary to the models of interest.
\end{enumerate}
\end{remark}

With this terminology in place, the main results itemized above are reflected in Table \ref{T:table}.  (The second and fourth rows should be read right-to-left to match the sense of the corresponding results.)

In Section \ref{S:logic} we indicate how pairs of properties can be identified by the use of Henson's \textit{positive bounded logic}, a version of model theory which interacts well with the ultraproducts of functional analysis.  At this writing more groundwork is required before positive bounded logic, or one of its equivalents, can be applied to the situations of this paper.  We expect the interaction between operator algebras and model theory to be fruitful in the near future.


\begin{table} \label{T:table}
\begin{center}
\begin{tabular}{|c|c|}
\hline
\textbf{Model version} & \textbf{Ultrapower version} \\
\textbf{(``approximative")} & \textbf{(``precise")} \\ \hline
approximate innerness & $\aleph_0$-local innerness \\ \hline
innerness & innerness \\ \hline
weak approximate unitary equivalence & unitary equivalence \\ \hline
``approximate pointwise innerness" (universal) & pointwise innerness \\ \hline
\end{tabular} \caption{Properties of models, with corresponding properties of their ultrapowers, in the context of this paper.}
\end{center}
\end{table}

\smallskip

\textbf{Acknowledgments.} We thank Nate Brown, Ward Henson, Narutaka Ozawa, and Nik Weaver for valuable comments.


\section{Relations between an automorphism and its ultrapower} \label{S:aut}

A family $\{\theta_j\} \subset \text{Aut}(\M)$ determines an automorphism $\Pi \theta_j$ of $\ell^\infty(\M)$ which descends to a well-defined automorphism of $\M^\omega$,
$$(\theta_j): (x_j) \mapsto (\theta_j(x_j)).$$
This automorphism is called an ultraproduct of the $\{\theta_j\}$, and one can only recover the representing sequence $(\theta_j)$ up to an obvious equivalence relation.  (The term ``liftable," which would seem appropriate here, has a different established meaning for automorphisms of ultrapowers (\cite[Section 5.2]{O}).)  It is natural to wonder whether every automorphism of $\M^\omega$ is such an ultraproduct.  We will show elsewhere, in joint work with Ilijas Farah, that sometimes the answer is negative.  (The first version of this result was proved by Farah and Nik Weaver.)


In any case, the subgroup $\{(\theta_j) \mid \{\theta_j\} \subset \text{Aut}(\M)\} < \text{Aut}(\M^\omega)$ itself contains two distinguished subgroups: those for which each $\theta_j$ is some fixed $\theta$ (in which case we denote the ultrapower automorphism by $\theta^\omega$), and those for which each $\theta_j$ is inner.  The latter is nothing but the inner automorphisms of $\M^\omega$.  We determine the intersection of these two subgroups in Theorem \ref{T:inner} below.

For countably-generated factors, the approximate innerness of $\theta$ amounts to the fact that $\theta^\omega$ agrees with an inner automorphism of $\M^\omega$ on the subalgebra $\M$.  This is well-known, but here we break up the logic in order to emphasize the roles of cardinality and local innerness.

\begin{theorem} \label{T:ai}
For an automorphism $\theta$ of a $\text{II}_1$ factor $\M$, these conditions are equivalent:
\begin{enumerate}
\item $\theta$ is approximately inner;
\item $\theta^\omega$ is $\aleph_0$-locally inner;
\item $\theta^\omega$ is approximately inner.
\end{enumerate}
The following condition implies the previous ones and is equivalent to them if $\M$ is countably-generated, but not in general:
\begin{enumerate}
\item[(4)] $\theta^\omega$ agrees with an inner automorphism of $\M^\omega$ on $\M$.

\end{enumerate}
\end{theorem}

\begin{proof} (1) $\Rightarrow$ (2): Given a countable family $\{X^n\} = \{(x^n_j)\}$, use the approximate innerness of $\theta$ to find unitaries $u_j$ such that
$$\|\theta(x^n_j) - u_j x^n_j u_j^*\|_2 \leq 2^{-j}, \qquad \forall n \leq j.$$
Then for each $n$,
$$\text{Ad}((u_j))(X^n) = (u_j x^n_j u_j^*) = (\theta(x^n_j)) = \theta^\omega(X^n).$$

(3) $\Rightarrow$ (1): Given $\varepsilon > 0$ and a finite set $\{x^n\} \subset \M \subset \M^\omega$, by the approximate innerness of $\theta^\omega$ we can find $U = (u_j) \in \U(\M^\omega)$ such that
$$\frac{\varepsilon}{2} \geq \max_n \|\theta^\omega((x^n)) - U(x^n)U^*\|_2 = \max_n \lim_{j \to \omega} \|\theta(x^n) - u_j x^n u_j^*\|_2.$$
Thus there must be some index $j_0$ with
$$\max_n \|\theta(x^n) - u_{j_0} x^n u_{j_0}^*\|_2 < \varepsilon,$$
as required.

(4) $\Rightarrow$ (1): Given $\varepsilon > 0$ and a finite set $\{x^n\} \subset \M \subset \M^\omega$, use the hypothesis to find $(u_j) \in \U(\M^\omega)$ such that $\theta^\omega$ agrees with $\text{Ad}((u_j))$ on all $(x^n)$.  This means
$$0 = \lim_{j \to \omega} \|\theta(x^n) - u_j x^n u_j^*\|_2, \qquad \forall n.$$
Again there must be some index $j_0$ with
$$\max_n \|\theta(x^n) - u_{j_0} x^n u_{j_0}^*\|_2 < \varepsilon,$$
as required.

The implications [(2) $\Rightarrow$ (3)] and (under the hypothesis that $\M$ is countably-generated) [(2) $\Rightarrow$ (4)] are trivial.

In \cite[Theorem 2.5]{S2008}, we displayed an outer $\aleph_0$-locally inner automorphism $\theta$ of a $\text{II}_1$ factor $\M$ which was constructed as a union $\cup_{\alpha < \aleph_1} \M_\alpha$.  This satisfies (1), since $\aleph_0$-local innerness implies approximate innerness.  We indicate why (4) fails, referring the reader to \cite{S2008} for supporting details.  It is required to show that an automorphism of the form $\text{Ad}((u_j))$ cannot agree with $\theta^\omega$ on all of $\M \subset \M^\omega$.  Let each $u_j \in \M_{\alpha_j}$, and set $\beta = (\sup \alpha_j) + 1 < \aleph_1$.  By construction $\M_\beta = \mathbb{M}_2 \otimes \M_{(\sup \alpha_j)}$; let $x \in \M_\beta$ be the element $(\begin{smallmatrix} 1 & 0 \\ 0 & 0 \end{smallmatrix}) \otimes 1$.  Again by construction $\theta(x) \neq x$.  Thus, viewing $x \in \M \subset \M^\omega$,
$$\text{Ad}((u_j))((x)) = (u_j x u_j^*) = (x) \neq (\theta(x)) = \theta^\omega((x)). \qedhere$$
\end{proof}



\begin{remark} \label{T:first} ${}$

\begin{enumerate}
\item The equivalence of (2) and (3) in fact holds for all automorphisms of $\M^\omega$ (Corollary \ref{T:precise}).

\item Haagerup and St\o rmer showed that if $\theta^\omega$ is inner, then $\theta$ is approximately inner (\cite[Theorem 6.2]{HSpinn}).  Theorem \ref{T:ai} gets the same conclusion from a weaker hypothesis, approximate innerness of $\theta^\omega$.   This suggests that their stronger hypothesis should have a stronger conclusion, and indeed it does (Theorem \ref{T:inner}).

\item Since the implication [(1) $\Rightarrow$ (2)] is based on an ultrafilter of $\mathbb{N}$, one might hope that for a suitable ultrafilter of a larger index set, the ultrapower of an approximately inner automorphism may actually be inner.  This does not happen (Theorem \ref{T:inner}).
\end{enumerate}
\end{remark}


Some parts of the next lemma are probably known, but we lack a reference.   We thank Narutaka Ozawa for suggesting the use of Dixmier averaging, which simplified our original argument.

We follow the convention of denoting the norm in $\B(\M, L^2(\M))$ by $\|\cdot\|_{\infty, 2}$.

\begin{lemma} \label{T:gp}
Let $\M$ be a $\text{II}_1$ factor.
\begin{enumerate}
\item For any $u \in \U(\M)$ we have
$$\|u - \tau(u)1 \|_2 \leq \|\textnormal{\mbox{Ad}}(u) - \textnormal{\mbox{id}}\|_{\infty, 2}.$$
\item Consider the following groups equipped with metrics: $($inner automorphisms of $\M$, $\|\cdot\|_{\infty, 2})$ and $((\U(\M)/\T)$ = projective unitary group of $\M$, quotient of the $L^2$ metric$)$.  Writing $\overline{u}$ for the coset of $u$ in $(\U(\M)/\T)$, the group isomorphism $\overline{u} \leftrightarrow \textnormal{\mbox{Ad}}(u)$ is Lipschitz continuous.
\item The inner automorphisms of $\M$ are complete in $\|\cdot\|_{\infty, 2}$.
\end{enumerate}

\end{lemma}

\begin{proof}
Dixmier's averaging theorem (\cite[Th\'{e}or\`{e}me 12]{D1949}) says that $\overline{\text{conv}\{vuv^* \mid v \in \U(\M)\}}^{\|\|} \cap \C = \{\tau(u)1\}$.  So for any $\varepsilon > 0$ we can find unitaries $v_1, \dots v_n$ and positive scalars $c_1, \dots c_n$ with $\sum c_n = 1$ such that $\|\sum c_j v_j u v_j^* - \tau(u)1\|_2 < \varepsilon$.  (We only require the $L^2$ estimate.)  Then we have the $L^2$ approximations
$$\tau(u)1 \overset{\varepsilon}{\approx} \sum c_j v_j u v_j^* = \sum c_j v_j (u v_j^* u^*) u \overset{\|\text{Ad}(u) - \text{id}\|_{\infty, 2}}{\approx} \sum c_j v_j (v_j^*) u = u.$$
Since this is true for any $\varepsilon$, we obtain (1).

The function $\text{dist}(\overline{u}, \overline{v}) = \min_{|\lambda|=1}\|u - \lambda v\|_2$ defines a metric on $(\U(\M)/\T)$ (essentially because $\T$ is a closed subgroup whose multiplicative actions on $\U(\M)$ are isometric).  We compute
\begin{align} \label{E:lip1}
[\text{dist}(\overline{u}, \overline{v})]^2 &= \min_{|\lambda|=1} \|u - \lambda v\|_2^2 \\ &= \notag \min_{|\lambda|=1} \|uv^* - \lambda\|_2^2 \\ &= \notag \min_{|\lambda|=1} \|(uv^* - \tau(uv^*)) + (\tau(uv^*) - \lambda)\|_2^2 \\ &= \notag \|uv^* - \tau(uv^*)\|^2 + \min_{|\lambda|=1}\|\tau(uv^*) - \lambda\|_2^2 \\ &= \notag \|uv^* - \tau(uv^*)\|^2 + (1 - |\tau(uv^*)|)^2 \\ &\leq \notag \|uv^* - \tau(uv^*)\|^2 + (1 - |\tau(uv^*)|^2) \\ &= \notag 2\|uv^* - \tau(uv^*)\|^2 \\ &\leq \notag 2 \|\text{Ad}(uv^*) - \text{id}\|_{\infty, 2}^2 \\ &= \notag 2 \|\text{Ad}(u) - \text{Ad}(v)\|_{\infty, 2}^2.
\end{align}
The fourth step is justified because $(uv^* - \tau(uv^*))$ is orthogonal to the scalars, and the second-to-last step is (1).
By an easy use of the triangle inequality, $\|\text{Ad}(u) - \text{Ad}(v)\|_{\infty, 2} \leq 2 \|u - v\|_2$, and this remains true if $v$ is multiplied by any unit scalar:
\begin{equation} \label{E:lip2}
\|\text{Ad}(u) - \text{Ad}(v)\|_{\infty, 2} \leq \min_{|\lambda| = 1} 2 \|u - \lambda v\|_2 = 2 \text{ dist}(\overline{u}, \overline{v}).
\end{equation}
Lipschitz continuity of the group isomorphism in (2) follows from \eqref{E:lip1} and \eqref{E:lip2}.


For (3), it suffices by (2) to show that $(\U(\M)/\T)$ is complete in the quotient of the $L^2$ metric.  So let $\{\overline{u_j}\} \subset (\U(\M)/\T)$ be a Cauchy sequence.  Choose a subsequence with $\text{dist}(\overline{u_{j_k}}, \overline{u_{j_{k-1}}}) < 2^{-k}$.  Multiplying each $u_{j_k}$ in turn by an appropriate unimodular scalar (and still denoting the sequence by $\{u_{j_k}\}$), we may obtain $\|u_{j_k} - u_{j_{k-1}}\|_2 < 2^{-k}$. By $L^2$ completeness of $\U(\M)$, the sequence $\{u_{j_k}\}$ converges in $L^2$ to some $u \in \U(\M)$.  It is immediate that $\{\overline{u_{j_k}}\}$ converges to $\overline{u}$ in the quotient metric, and the Cauchy sequence $\{\overline{u_j}\}$ must converge to $\overline{u}$ as well.  (A general fact, proved identically and for which we have no reference: when a group is a complete metric space, and a closed subgroup acts isometrically by right multiplication, then the left coset space is complete in the quotient metric.)
\end{proof}

\begin{remark}
The constant 2 in inequality \eqref{E:lip2} is sharp, and the isomorphism in Lemma \ref{T:gp}(2) is not isometric.  This can be verified with unitaries from a copy of $\mathbb{M}_2$ inside $\M$.  The inequalities in Lemma \ref{T:gp}(1) and \eqref{E:lip1} are probably true with better constants.
\end{remark}

\begin{theorem} \label{T:inner}
Let $\theta$ be an automorphism of the $\text{II}_1$ factor $\M$.  If $\theta^\omega$ is inner, then $\theta$ is inner.
\end{theorem}

\begin{proof}
Suppose $\theta^\omega = \text{Ad}((u_j))$.  For each $j$ find a contraction $x_j \in \M$ such that
$$\|[\theta - \text{Ad}(u_j)](x_j)\|_2 \geq \frac12 \|\theta - \text{Ad}(u_j)\|_{\infty, 2}$$
and compute
$$0 = \|\theta^\omega((x_j)) - \text{Ad}((u_j))((x_j))\|_2 = \lim_{j \to \omega} \|[\theta - \text{Ad}(u_j)](x_j)\|_2 \geq \frac12 \lim_{j \to \omega} \|\theta - \text{Ad}(u_j)\|_{\infty, 2}.$$
Lemma \ref{T:gp}(3) then implies that $\theta$ is inner.
\end{proof}



Theorem \ref{T:inner} is valid for any ultrafilter, not just $\omega \in (\beta \mathbb{N} \setminus \mathbb{N})$, corroborating Remark \ref{T:first}(3).

From Theorem \ref{T:ai} and Theorem \ref{T:inner}, we obtain new examples of outer (even $\aleph_0$-)locally inner automorphisms: $\theta^\omega$, for $\theta$ outer and approximately inner.  As mentioned earlier, any countably-generated $\Gamma$ factor admits such $\theta$.




\section{Approximate equivalence for maps into ultrapowers} \label{S:ae}

We first specify some notation.  For $y$ in a von Neumann algebra, we write $s_\ell(y)$ for the \textit{left support} of $y$, which is the least projection $p$ with $py = y$.  When the algebra is represented on a Hilbert space the left support is nothing but the range projection.

For unital *-homomorphisms $\pi, \rho$ from a C*-algebra $\A$ to a von Neumann algebra $\N$, we consider the following four relations.

\begin{enumerate}
\item \textit{Unitary equivalence}: $\exists u \in \U(\N), \: (\text{Ad }u ) \circ \pi = \rho$.
\item \textit{Approximate unitary equivalence}: $\exists \{u_\alpha\} \subset \U(\N), \: (\text{Ad } u_\alpha) \circ \pi \to \rho$ in the point-norm topology.
\item \textit{Weak approximate unitary equivalence}: $\exists \{u_\alpha\} \subset \U(\N), \: (\text{Ad }u_\alpha ) \circ \pi \to \rho$, and $\exists \{v_\alpha\} \subset \U(\N), \: (\text{Ad }v_\alpha ) \circ \rho \to \pi$, both in the point-weak topology.  (It makes no difference to use the point-strong or point strong*, as first pointed out in \cite[Section 1]{H1981}.)
\item \textit{Equal rank}: for all $x \in \A$, $s_\ell(\pi(x)) \sim s_\ell(\rho(x))$.
\end{enumerate}

Obviously conditions (1) through (3) are progressively weaker.  We always have the implication [(2) $\Rightarrow$ (4)], but [(3) $\Rightarrow$ (4)] holds if and only if $\N$ is a direct sum of $\sigma$-finite von Neumann algebras.  To see these implications, note that $s_\ell(\pi(x)) = \chi_{\C \setminus \{0\}}(\pi(x^*x))$ and apply \cite[Theorem 5.4]{S2007b}.  As for the failure of the second implication when $\N$ is not a direct sum of $\sigma$-finite algebras, consider $\pi,\rho: \C^2 \to \N$ such that $\pi(1 \oplus 0)$ and $\rho(1 \oplus 0)$ have equal central support, while one is $\aleph_1$-homogeneous and the other is $\aleph_0$-homogeneous (\cite[Theorem 3.5 and Proposition 3.8]{S2007a}).


There is much left to understand about the partial validity of the converses to these implications.  We think of these as generalizations of Voiculescu's noncommutative Weyl-von Neumann theorem (\cite[Theorem 1.5]{V}), as Hadwin's beautiful reformulation (\cite[Theorem 3.14]{H1981}) says that [(4) $\Rightarrow$ (2)] when $\N = \B(\h)$.  In general the implication [(4) $\Rightarrow$ (3)] can fail even for $\A$ separable and $\N = \RR$ (\cite[Corollary 3.5]{H1998}).   See \cite{H1998, DH, S2007b} for more discussion.

In case $\N$ is a $\text{II}_1$ factor, the conditions above can be simplified and related to other familiar terms.  The net $\{v_\alpha\}$ is not needed in (3), as
$$\|u_\alpha \pi(x) u_\alpha^* - \rho(x)\|_2 \to 0 \: \iff \: \|\pi(x) -  u_\alpha^* \rho(x) u_\alpha\|_2 \to 0.$$
Also (4) is the same as requiring $\tau \circ \pi = \tau \circ \rho$ (\cite[Lemma 3]{DH}).  Viewing an automorphism as a *-homomorphism from $\N$ to $\N$, approximate innerness amounts to being weakly approximately unitarily equivalent to the identity.  And any pair of automorphisms satisfy (4), by uniqueness of the trace.  Actually the only factors which admit ``rank-changing" automorphisms are those $\text{II}_\infty$ factors whose fundamental group is nontrivial.

In the rest of this section we focus on $\N = \M^\omega$.  The main results are that [(3) $\Rightarrow$ (1)] for $\A$ separable and, using a result of Ding and Hadwin, [(4) $\Rightarrow$ (1)] if in addition $\A$ is nuclear.
\begin{theorem} \label{T:ue}
Let $\A$ be a separable C*-algebra, $\M$ be a $\text{II}_1$ factor, and $\pi,\rho: \A \to \M^\omega$ be unital *-homomorphisms.  If $\pi$ and $\rho$ are weakly approximately unitarily equivalent, then they are unitarily equivalent.
\end{theorem}

\begin{proof}
Let $\{x^k\} \subset \A_1$ be a countable generating set for $\A$, and let $\pi(x^k) = A^k = (a^k_j)$, $\rho(x^k) = B^k = (b^k_j)$.  Although nets are unavoidable in the general definition of weak approximate unitary equivalence, here a sequence will work.  For each $n \in \mathbb{N}$ find a unitary $U^n$ with
$$\max_{1 \leq k \leq n} \|U^n A^k U^{n*} - B^k\|_2 \leq \frac1n.$$
It is straightforward to check that $\text{Ad }(U^n)\circ \pi \to \rho$ in the point-strong topology.

Now let $U^n = (u^n_j)$.  For $k \leq n$ we have
$$\lim_{j \to \omega} \|u^n_j a^k_j u^{n*}_j - b^k_j\|_2 \leq \frac1n.$$

We consider the function
$$f: \U(\M^\omega) \to [0,2],$$
$$W = (w_j) \mapsto \sum_k 2^{-k} \|W A^k W^* - B^k\|_2 = \lim_{j \to \omega} \sum_k 2^{-k} \|w_j a^k_j w^*_j - b^k_j\|_2.$$
(As for the interchange of limits, just note that the functions $\{(j \mapsto 2^{-k} \|w_j a^k_j w^*_j - b^k_j\|_2)\}_k$ are absolutely summable in $\ell^\infty = C(\beta \mathbb{N})$, and we are evaluating them at $j = \omega$.)  By construction
$$f(U^n) = \sum_k 2^{-k} \|U^n A^k U^{n*} - B^k\|_2 \leq \left(\sum_{k=1}^n 2^{-k} \frac1n \right) + \left(\sum_{k = n+1}^\infty 2^{-k} \cdot 2 \right) \leq \frac1n + 2^{-n+1}.$$
So we have for each $n$,
$$\lim_{j \to \omega} \left(\sum_k 2^{-k} \|u^n_j a^k_j u^{n*}_j - b^k_j\|_2 \right) \leq \frac1n + 2^{-n+1}.$$
For each $j$, define $v_j$ from among $u^1_j, u^2_j, \dots, u^j_j$ so that the quantity in parentheses above is minimized.  Thus for $j \geq n$,
$$\left(\sum_k 2^{-k} \|v_j a^k_j v^*_j - b^k_j\|_2 \right) \leq \left(\sum_k 2^{-k} \|u^n_j a^k_j u^{n*}_j - b^k_j\|_2 \right).$$
Taking limits, we conclude that for any $n$,
$$\lim_{j \to \omega} \left(\sum_k 2^{-k} \|v_j a^k_j v^*_j - b^k_j\|_2 \right) \leq \frac1n + 2^{-n+1}.$$
Thus with $V = (v_j) \in \M^\omega$, we must have $f(V) = 0$.  But then $V A^k V^* = B^k$ for all $k$, so that $\pi$ and $\rho$ are unitarily equivalent.
\end{proof}


\begin{corollary} \label{T:precise}
Let $\M$ be a $\text{II}_1$ factor.  If an automorphism of $\M^\omega$ is approximately inner, then it is $\aleph_0$-locally inner.
\end{corollary}

\begin{proof}
Let $\alpha$ be an approximately inner automorphism of $\M^\omega$.  Take a countable family $\{X^j\} \subset \M^\omega$, and set $\A = C^*(\{X^j\})$.  Then the two representations $\text{id}, \alpha \circ \text{id}: \A \to \M^\omega$ are weakly approximately unitarily equivalent.  By Theorem \ref{T:ue} they are actually unitarily equivalent, so that $\alpha$ agrees with some inner automorphism on all the $X^j$.
\end{proof}

\begin{theorem} \label{T:dh} $($\cite[Theorem 5]{DH} or \cite[Theorem 2.1]{H1998}$)$
Let $\A$ be a nuclear C*-algebra, $\M$ be a $\text{II}_1$ factor, and $\pi,\rho: \A \to \M$ be unital *-homomorphisms with $\tau \circ \pi = \tau \circ \rho$.  Then $\pi$ and $\rho$ are weakly approximately unitarily equivalent.
\end{theorem}
(The theorem is stated in \cite{DH} with the assumption that $\M$ acts on a separable Hilbert space, but this is only so that disintegration theory may be applied to a non-factor.)

\begin{corollary} \label{T:nuc}
Let $\A$ be either a separable nuclear C*-algebra or a countably-generated hyperfinite finite von Neumann algebra, $\M$ be a $\text{II}_1$ factor, and $\pi,\rho: \A \to \M^\omega$ be unital (normal, in the von Neumann algebra case) *-homomorphisms with $\tau_{\M^\omega} \circ \pi = \tau_{\M^\omega} \circ \rho$.  Then $\pi$ and $\rho$ are unitarily equivalent.
\end{corollary}

\begin{proof}
For the C*-algebra version, combine Theorems \ref{T:dh} and \ref{T:ue}.  For the von Neumann algebra version, just create a weakly dense nuclear C*-subalgebra as the norm closure of a weakly dense increasing sequence of finite-dimensional C*-subalgebras.
\end{proof}

\begin{remark}
Let $\N$ be a finitely-generated von Neumann algebra which embeds in $\RR^\omega$, and fix a faithful normal trace $\tau_\N$ on $\N$.  Jung proved in \cite{J} that $\N$ is hyperfinite if and only if all *-homomorphisms $\pi: \N \to \RR^\omega$ satisfying $\tau_{\RR^\omega} \circ \pi = \tau_\N$ are unitarily equivalent.  The forward implication of Jung's theorem was also established in \cite[Theorem 6.1]{FGL} without the hypothesis of finite generation.  These results contain Corollary \ref{T:nuc} for the situation $\M = \RR$.

Note that the C*-algebra version of Corollary \ref{T:nuc} also follows easily from the von Neumann algebra version, as $\overline{\pi(\A)}^s \simeq \overline{\rho(\A)}^s$ is a hyperfinite von Neumann algebra to which $\pi$ and $\rho$ extend (because $\tau \circ \pi = \tau \circ \rho$).
\end{remark}

Haagerup and St\o rmer showed in \cite[Theorem 6.2]{HSpinn} that automorphisms of $\M^\omega$ of the form $\theta^\omega$ are always pointwise inner.  The first statement of Corollary \ref{T:pinn} shows that all automorphisms of $\M^\omega$ actually have a stronger property.



\begin{corollary} \label{T:pinn}
Let $\M$ be a $\text{II}_1$ factor.
\begin{enumerate}
\item Let $\A$ be a separable nuclear C*-subalgebra of $\M^\omega$, and let $\alpha$ be an automorphism of $\M^\omega$.  Then $\alpha$ agrees with some inner automorphism on $\A$.
\item Every automorphism of $\M^\omega$ is pointwise inner.
\item Let $\N$ be a countably-generated hyperfinite von Neumann subalgebra of $\M^\omega$, and let $\alpha$ be an automorphism of $\M^\omega$.  Then $\alpha$ agrees with some inner automorphism on $\N$.
\end{enumerate}
\end{corollary}

\begin{proof}
For the first part, unitize $\A$ if necessary and then note that $\text{id}, \alpha \circ \text{id}: \A \to \M^\omega$ satisfy the conditions of Corollary \ref{T:nuc}.  The second part is an immediate consequence of the first.  As in the proof of Corollary \ref{T:nuc}, the third part follows from (also implies) the first.
\end{proof}

\begin{remark}
The second statement essentially follows from the proof of \cite[Lemma 7.1]{P1983} or \cite[Lemma 4.2]{FGL}, although there only diffuse abelian subalgebras are discussed.  Note that the statement becomes false if ``pointwise inner" is replaced with ``locally inner" (\cite[Proposition 2.2]{S2008}).

The third statement is false without the hyperfiniteness assumption, as the ultrapower of a non-approximately inner automorphism of $\M$ fails to be $\aleph_0$-locally inner (Theorem \ref{T:ai}).  Since countably-generated hyperfinite von Neumann algebras are singly-generated (\cite[Theorem 1]{SS}), one may view the asserted property of $\alpha$ as somewhere between pointwise innerness and local innerness.  In the taxonomy of \eqref{E:impl}, this is ``hyperfinite-local innerness", where $1 > \text{``hyperfinite"} > \frac12$.
\end{remark}

\section{A logical conclusion} \label{S:logic}

There are a few different versions of model theory which interact well with the ultrapowers of functional analysis.  In this section we give a simplified idea of Henson's positive bounded logic -- see \cite{HI} for the whole story -- then sketch how ``approximative/precise" property pairs can be justified by metamathematical theorems.  At this writing tracial von Neumann algebras have not been treated in the published literature as a class of model structures, although there is current work undertaking to do so using the so-called ``model theory for metric structures" (\cite{BBHU}), which may be viewed as a generalization of positive bounded logic.  See \cite{BHJR} for the assertion (without proof) that tracial von Neumann algebras are an axiomatizable class.

As in classical model theory, one starts with a language which is suitable for describing the functions and relations of the models one is interested in, e.g. addition, scalar multiplication, and norm functions for Banach spaces.  Syntactically one is limited to the \textit{positive bounded formulas}: the ones which can be built out of non-strict norm inequalities via $\wedge$, $\vee$, and quantification over bounded sets -- implication and negation are off limits.  Given a model $M$ and a positive bounded sentence (=formula with no free variables) $\varphi$, we write $M \models \varphi$ and say that $M$ \textit{satisfies} $\varphi$ if $\varphi$ is true in $M$.  Now suppose only that $M$ satisfies all sentences obtained by weakening the constants of $\varphi$ by arbitrarily small amounts.  In this case we write $M \models_a \varphi$ and say that $M$ \textit{approximately satisfies} $\varphi$.

Positive bounded logic is a model theory in which approximate satisfaction is used in place of satisfaction.  There are analogues of many of the classical theorems: compactness, L\"{o}wenheim-Skolem, \L o\'s, Keisler-Shelah, etc.  Regarding ultrapowers, one can deduce (\cite[Corollary 9.3 and Proposition 9.26]{HI})
\begin{equation} \label{E:updown}
M \models_a \varphi \quad \iff \quad M^\omega \models \varphi \quad \iff \quad M^\omega \models_a \varphi.
\end{equation}
This is a version of our pattern (but not the whole story).

Below we apply positive bounded logic to revisit two examples based on $\text{II}_1$ factors equipped with the $L^2$ norm.  Our discussion is conceptually accurate but brief, so we do not present the languages explicitly.  The second example is meant to reassure operator algebraists that they have already been working with approximate satisfaction for a long time.

\begin{example} \label{T:reproof} (Second proof of Theorem \ref{T:inner}) In the terminology of the Introduction, we need to show that the model version of innerness is again innerness.

Working in the class whose models consist of a $\text{II}_1$ factor with a single automorphism, the assertion that $(\M, \theta)$ is inner can be expressed as $(\M, \theta) \models \varphi$, where $\varphi$ is the positive bounded sentence
$$(\exists_1 u) (\forall_1 x) \|\theta(x) - uxu^*\|_2 \leq 0.$$
Here ``$\exists_1$" means ``there exists in the closed ball of radius 1"; ``$\forall_1$" is similar.  According to \eqref{E:updown} the model version of innerness is $(\M, \theta) \models_a \varphi$.  This means that $(\M, \theta)$ satisfies every sentence of the form
\begin{equation} \label{E:apsat}
(\exists_{1+\delta_1} u) (\forall_{1-\delta_2} x) \|\theta(x) - uxu^*\|_2 \leq \delta_3,
\end{equation}
for arbitrarily small $\delta_j$.  The $u$ whose existence is asserted by \eqref{E:apsat} need not be unitary, but taking scalar $x$ gives $uu^* \approx 1$ in $L^2$, which implies that there is a unitary near $u$ which will satisfy a slightly worse bound.  Then the sentences \eqref{E:apsat} say that $\theta$ is a uniform-$L^2$ limit of inner automorphisms.  The proof is completed by invoking Lemma \ref{T:gp}(3) as before.
\end{example}

\begin{example} \label{T:Gamma} (Property $\Gamma$)  Let $\M$ be a model for the class of $\text{II}_1$ factors, with language sufficiently rich to allow the positive bounded sentences
$$\varphi_n: \qquad (\forall_1 x_1) (\forall_1 x_2) \dots (\forall_1 x_n) (\exists_1 u) \left[ ( \|u^* u - 1\|_2 \leq 0 ) \wedge (|\tau(u)| \leq 0) \wedge \left( \bigwedge_{j=1}^n \|ux_j - x_j u\|_2 \leq 0 \right) \right].$$
The condition $\M \models \{\varphi_n\}_{n \in \mathbb{N}}$ says that the relative commutant of any finite set in $\M$ contains a unitary with zero trace; equivalently, the relative commutant of any finite set in $\M$ is nontrivial.  This is clearly impossible if $\M$ is finitely-generated.  In fact it is still impossible if $\M$ is countably-generated, as $\M$ contains an irreducible hyperfinite (so singly-generated) subfactor (\cite[Corollary 4.1]{P1981}).

The condition $\M \models_a \{\varphi_n\}$ says that $\M$ has property $\Gamma$.  (As in the previous example, one can $L^2$-perturb $u$ to a trace-zero unitary with slightly worse bounds on the commutators.)  By \eqref{E:updown} this is equivalent to $\M^\omega \models_a \{\varphi_n\}$, i.e. $\M$ has $\Gamma$ if and only if $\M^\omega$ has $\Gamma$ (\cite[Corollary 5.2]{FGL}).  Also by \eqref{E:updown}, it is equivalent to $\M^\omega \models \{\varphi_n\}$, which says that any finite subset of $\M^\omega$ has nontrivial relative commutant.  If $\M$ is finitely-generated, a small argument shows that this in turn is equivalent to $\M' \cap \M^\omega \neq \C$.  Well-known conclusion (for finitely-generated $\M$): $\M$ has $\Gamma$ if and only if $\M' \cap \M^\omega \neq \C$.

One can and should replace ``finite" with ``countable" in the last three sentences.  This requires a slightly different approach which appeals to the $\aleph_1$-saturation of $\M^\omega$, and the details are omitted.  Notice the similarity between the four conditions in the previous paragraph and those of Theorem \ref{T:ai}.  Also notice that the ultrapower version of property $\Gamma$ is ``countable subsets have nontrivial relative commutant."
\end{example}

Actually there are several places in the literature on operator algebraic ultrapowers where a logical approach would be effective, but in few, if any, of these cases would it be shorter -- especially considering the extra machinery which must be introduced to the reader.  Kirchberg simply proves the version of $\aleph_1$-saturation that he needs in an Appendix (\cite[Lemma A.1 and Remark A.2]{K}), although he does not name it as such.  However, in recent work with Farah, the author has used logical methods to obtain results on ultrapowers which may be (in some sense) inaccessible via analytic techniques.  These will appear at a later date.

\end{document}